\documentclass[prl,%
 reprint,
 amsmath,amssymb,
 aps,
]{revtex4-2}

\usepackage{graphicx}
\usepackage{dcolumn}
\usepackage{bm}

\usepackage{geometry}
\usepackage{array}
\usepackage{epsfig}
\usepackage{rotating}
\usepackage{amssymb}
\usepackage{subfigure}
\usepackage{makecell}
\usepackage{booktabs}
\usepackage{dsfont}
\usepackage{psfrag}
\usepackage{bbold}
\usepackage{xcolor}
\usepackage{multirow}
\usepackage{amsmath,euscript,array,mathrsfs}
\usepackage[utf8]{inputenc}
\usepackage[T1]{fontenc}
\usepackage{enumitem}
\usepackage[colorlinks=true,linktocpage=true,linkcolor=blue,citecolor=blue]{hyperref}
\usepackage{amsthm}
\usepackage{float}

\usepackage[section]{placeins}

\newcommand{\be}{\begin{equation}}
\newcommand{\ee}{\end{equation}}
\newcommand{\bea}{\begin{eqnarray}}
\newcommand{\eea}{\end{eqnarray}}

\begin{document}

\preprint{APS/123-QED}

\title{\textbf{Patterns in Growth and Distribution of Unbounded Prime Number Walks}}

\author{Alberto Fraile$^{*}$, Daniel Fernández$^{\dagger}$, Roberto Martínez$^{\ddagger}$, Theophanes E. Raptis$^{\mathsection}$}

\affiliation{$^*$Catalan Institute of Nanoscience and Nanotechnology (ICN2), BIST \& CSIC Campus UAB, Bellaterra, Barcelona 08193, Spain}
\affiliation{$^{\dagger}$Fléttan EO ehf., Environmental Monitoring Analysis Unit, Reykjavík, Iceland}
\affiliation{$^{\ddagger}$Universidad del País Vasco / Euskal Herriko Unibertsitatea, 48940. Bizkaia, Spain}
\affiliation{$^{\mathsection}$Information Physics Institute Network, Athens, Greece}

\email{$^*$albertofrailegarcia@gmail.com}
\email{$^\dagger$licadnium@gmail.com}
\email{$^\ddagger$rmartinezballarin@gmail.com}
\email{$^\mathsection$traptis@protonmail.com}

\date{\today}

\begin{abstract}
In our previous work, we defined a prime walk on a square grid and presented several intriguing numerical results. Here, we demonstrate the main conjecture presented there, namely, that the area covered by the prime walk is unbounded. Taking this fact into account, we examine in further detail the properties of the PW and explore new questions that arise naturally in this analysis.
\end{abstract}

\maketitle


\section{Introduction}
\label{sec:Intro}

The distribution of prime numbers has fascinated mathematicians since ancient times. The apparent interplay between randomness and hidden patterns has inspired continuous exploration, aiming to uncover deeper insights into their nature \cite{crandall2005prime, fine2016number}.

Random walks (RW) are ubiquitous in science, and important from both theoretical and practical perspectives \cite{weiss1982, weiss1994random, Redner_2001}. They are one of the most fundamental types of stochastic processes and can be used to model numerous phenomena \cite{klafter2011random, falco}. Besides, RWs have been studied for decades on both regular lattices and (especially in recent years) on networks with a variety of structures \cite{citeulike:218038, regnier23, lemke16, wolf14, perez14, southier23}.

In our previous work \cite{fraile2021primes}, we conducted a straightforward mathematical experiment, defining a “random” walk on a square grid. The walk was dictated by, and dependent on, the sequence of the last digits of the primes. The Prime Walk (PW) is defined by assigning \((x, y)\) positions to a sequence of positive integers, from 1 to \(N\), on the \(XY\)-plane.


It begins at \((0,0)\), assigned to \(N = 1\), and advances according to the last digit of the next number \(N + 1\), but only if it is prime. If not, the position remains unchanged. When \(N + 1\) is prime, the walker moves: down for digit 1, up for 3, right for 7, and left for 9.

Note that the last digits of prime numbers are always 1, 3, 7, or 9, with two exceptions: 2 and 5, which are ignored by the algorithm. The methodology is described in our previous paper \cite{fraile2021primes}. The code used to create the PW is provided in \cite{supp}.

Other choices are possible, but only two alternatives are intrinsically distinct, which we denote as A2 and A3. In compact notation, these alternatives are: (A1: \(1\downarrow,\,3\uparrow,\,7\rightarrow,\,9\leftarrow\)), (A2: \(1\rightarrow,\,3\uparrow,\,7\downarrow,\,9\leftarrow\)), and (A3: \(1\leftarrow,\,3\uparrow,\,7\rightarrow,\,9\downarrow\)). Observe that any other possible definition is equivalent, under rotations or reflections, to one of the previously mentioned definitions.

Therefore, combining in a single study RWs and prime numbers creates a novel platform for exploring the properties of the prime number sequence (mod 10), with potential implications in areas such as number theory, stochastic processes, and cryptography.

As will be shown, all three algorithms appear to behave similarly in the asymptotic limit. However, we emphasize that this observation is empirical rather than a rigorous guarantee valid for arbitrarily large \(N\). Indeed, discovering significant differences among these definitions under certain conditions would be intriguing.

In \cite{fraile2021primes}, we defined the area \(A(N)\) as the number of unique points visited by the PW after \(N\) steps, commonly referred to as the \emph{range} in the literature~\cite{vallois1996range}. This quantity has attracted considerable interest across diverse fields, including physics, chemistry, and ecology.

\begin{figure*}[t]
\includegraphics[width=0.9\textwidth]{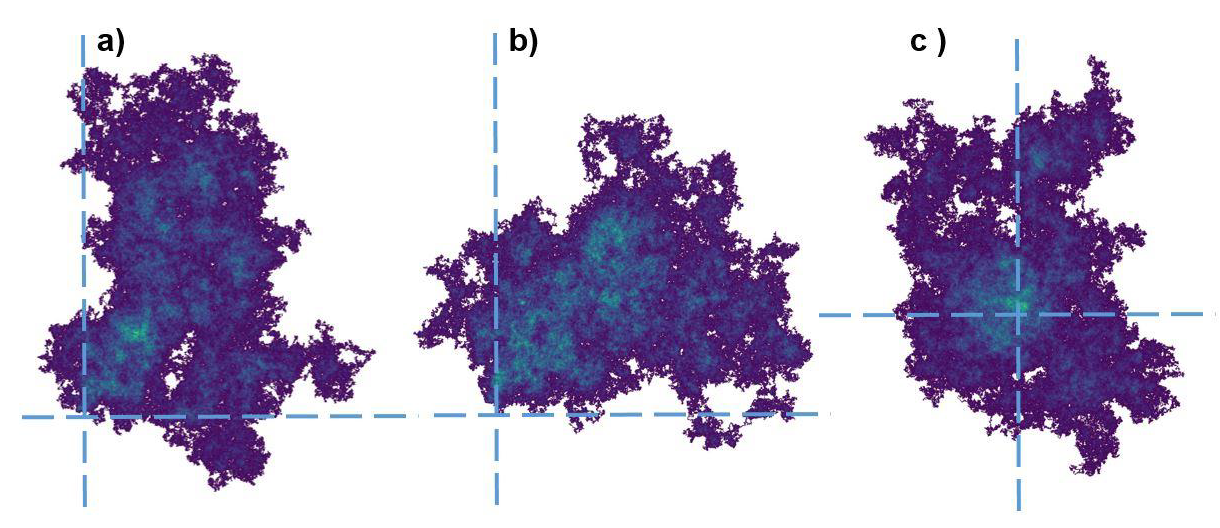}
\caption{\label{fig:area}Area covered by the three different prime walks: (a) A1, (b) A2, and (c) A3, up to \(N = 10^{10} \). Colors represent the number of visits, with brighter areas indicating higher visitation frequency. The dashed lines correspond to the coordinate axes.}
\end{figure*}

A classic theme in analytic number theory can be framed as exploring how the sequence of primes resembles (or deviates from) a random sequence of integers \cite{vineyard1963random, weiss1994random}. One of our key findings was that the PW covers an area approximately half that of a RW with the same number of steps, under certain conditions \cite{fraile2021primes}. This result is perplexing, and we were unable to identify a satisfactory explanation. Adding to the intrigue, over the explored interval (\(5 \times 10^{10}\)), the area covered is approximately one tenth the total number of primes in \(N\) steps - another unexpected result that remains unexplained.

Many conjectures emerge from these observations, and it is safe to say that most are extremely challenging to prove. In this Letter, we investigate further the properties of the PWs and, importantly, we prove the primary conjecture proposed in \cite{fraile2021primes}.

Figure~\ref{fig:area} illustrates the areas covered by the three distinct PWs, as defined previously, for \(N = 10^{10} \), corresponding to a total number of primes \(N_p = 455,052,509\). The color scale represents the frequency of visits to each site, with brighter colors indicating more frequently visited locations.

\section{Analysis of Area Growth}
\label{sec:AreaGrowth}

In our previous work~\cite{fraile2021primes}, we examined the PW up to \(N = 5 \times 10^{10}\) and observed that the area exhibited linear growth as a function of \(N\), where \(N\) represents the total length of the explored interval. An alternative analysis, however, can provide additional insights by plotting the results as a function of the number of steps, specifically the number of primes \(N_p\) within the interval \(N\). In either approach (the previous method or the one employed here), the computed area \(A(N)\) remains the same, since it is defined by the number of unique points \((x, y)\) visited.

\begin{figure}[H]
    \centering
    \includegraphics[width=0.5\textwidth]{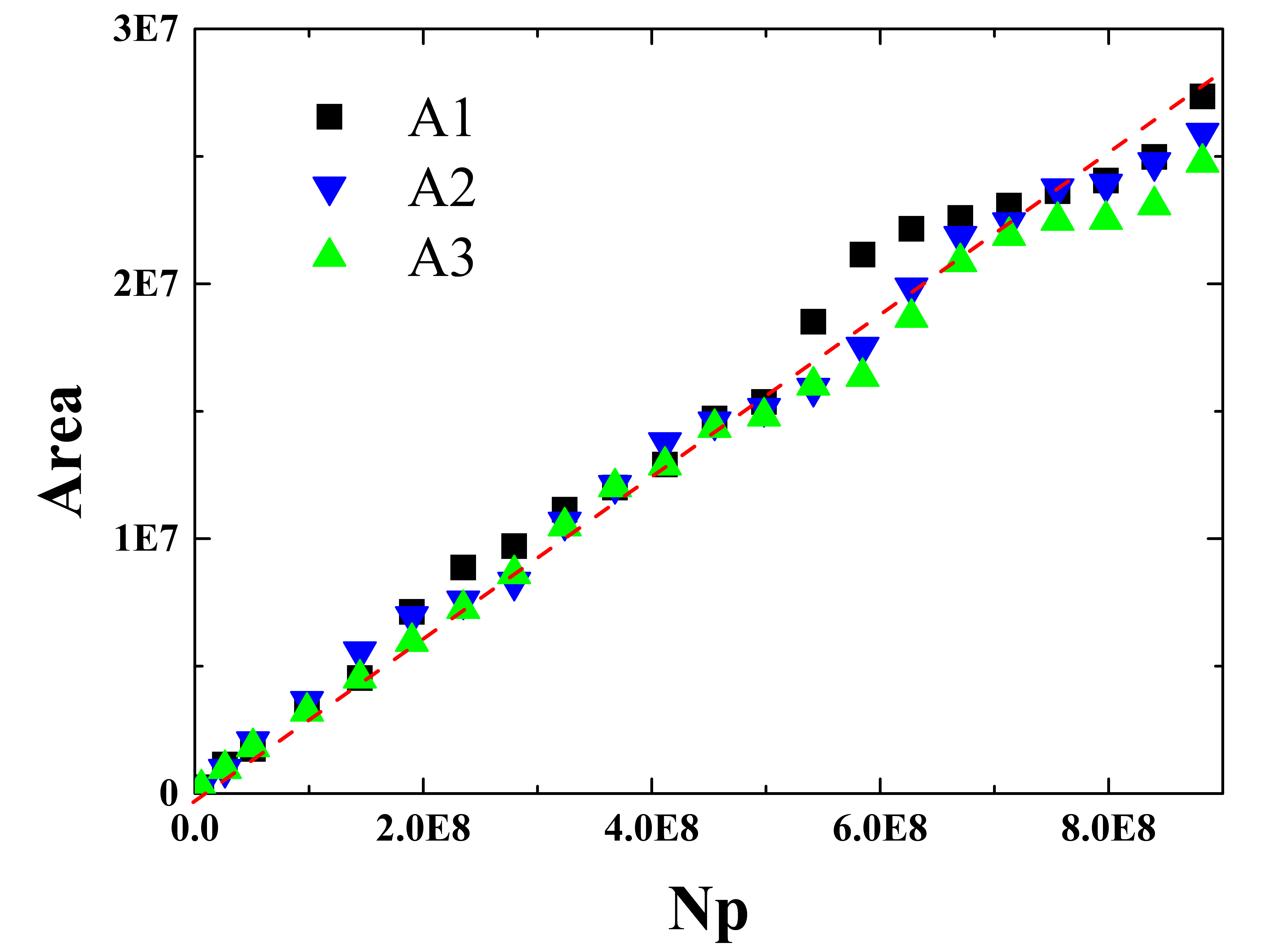}
    \caption{Area covered by the prime walks (A1: black squares, A2: blue inverted triangles, and A3: green triangles) versus the number of primes \(N_p\). The red line is included as a guide to the eye.}
    \label{fig:area_vs_primes}
\end{figure}

Figure~\ref{fig:area_vs_primes} shows the area covered by the three different prime walks (A1, A2, and A3) plotted against \(N_p\), up to \(N = 2\times 10^{10} \), corresponding to a total of \(N_p = 882,206,715\) primes.

Remarkably, despite the visibly different shapes of the visited regions (see Figure~\ref{fig:area}), the behavior of the three algorithms is strikingly similar, as demonstrated by their overlapping values. Performing a linear fit to each dataset yields nearly identical slopes: \(b_1 = 0.032 \pm 0.0004\), \(b_2 = 0.031 \pm 0.0003\), and \(b_3 = 0.029 \pm 0.0003\) (with asymptotic standard errors of approximately \(1.4\%\), \(1.03\%\), and \(1.097\%\), respectively). This strongly suggests that, in the asymptotic limit, the slopes might converge toward a common value, denoted here by the average \(\beta \approx 0.03\). Whether this value relates to a fundamental mathematical constant or represents a deeper underlying principle remains an intriguing open question.

\section{Proof of Infinite Growth}
\label{sec:Conjectures}

In \cite{fraile2021primes}, we conjectured that the area \(A(N)\) grows indefinitely as \(N\) increases. Since we can assume that the PW takes approximately the same number of steps in every direction, it seems reasonable to expect that the PW will hover around the origin without straying too far. However, this does not guarantee that the area is bounded. It is simple to construct cases in which the area grows indefinitely, even with an equal number of moves in all directions. For example, a spiral-like growth would result in a continuously expanding area, with the number of moves in each direction roughly equal for every turn.

One way to approach this problem is by examining sequences of consecutive primes ending in 1, 3, 7, or 9 within a given interval. The longest such sequence within an interval would determine the maximum possible displacement in a particular direction. If arbitrarily long sequences of primes with identical ending digits exist, the area \(A(N)\) must be unbounded, as these sequences would allow arbitrarily large steps in one or more directions.

\begin{figure}
    \centering
    \includegraphics[width=0.5\textwidth]{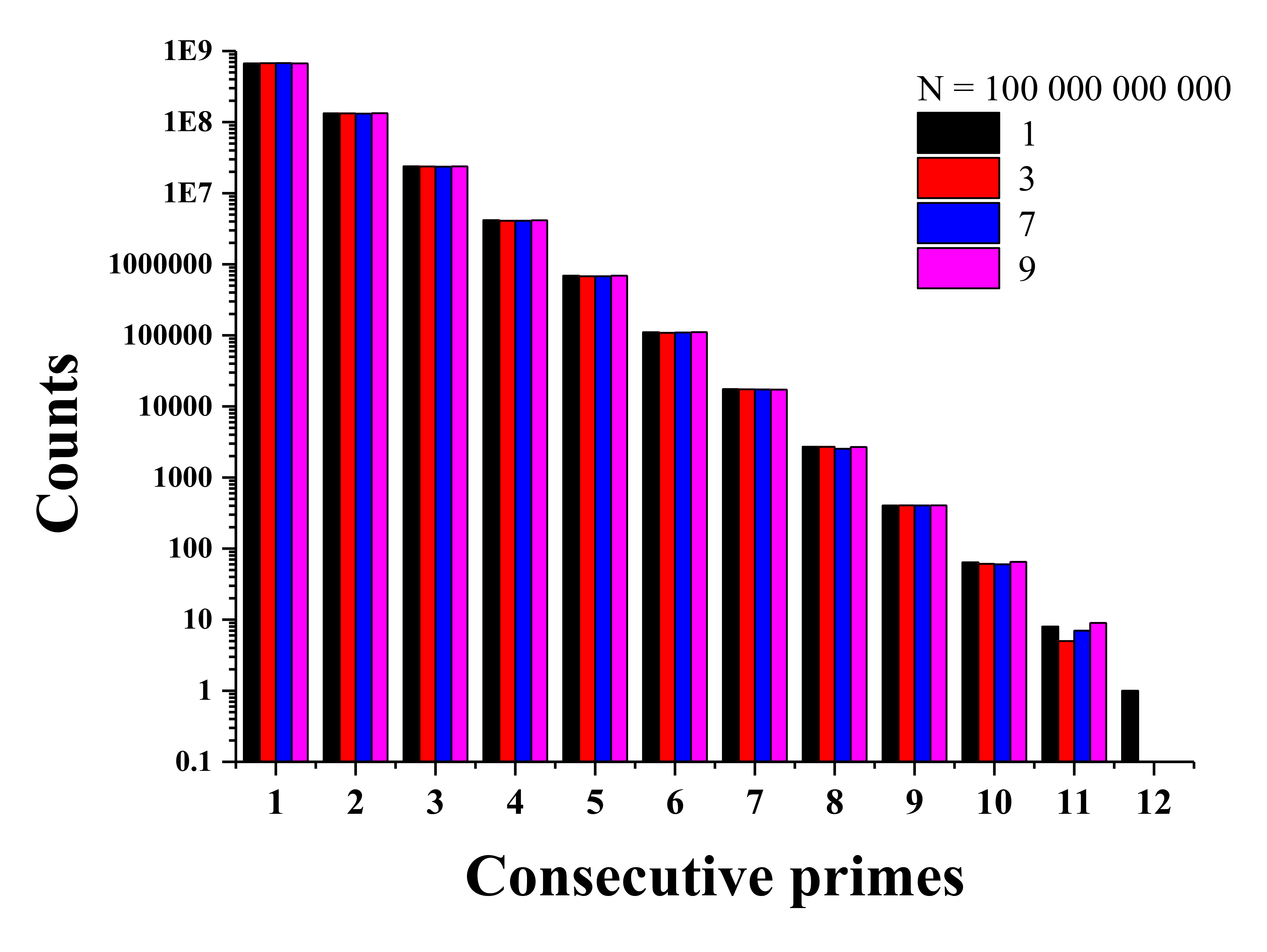}
    \caption{Number of consecutive primes ending in 1, 3, 7, and 9 within \(N = 10^{11}\). The maximum sequence (12 consecutive primes) occurs only once for primes ending in 1.}
    \label{fig:consecutive_primes}
\end{figure}

Indeed, numerical evidence supports this idea. For example, within the interval \(N = 10^{11}\), the longest sequence of consecutive primes ending in the same digit (1 in this case) is 12, and this occurs only once. No such sequences of length 12 are observed for primes ending in 3, 7, or 9. In other words, long sequences are rare, but it appears that such sequences inevitably arise as the interval becomes sufficiently large. The numerical data used to create Figure~\ref{fig:consecutive_primes} is provided as Table 1 in \cite{supp}. 

On the other hand, short sequences (those consisting of only 1 or 2 consecutive primes sharing the same last digit) overwhelmingly dominate, accounting for more than 95\% of all observed cases.

Upon examining Figure~\ref{fig:consecutive_primes}, it appears reasonable to speculate that up to \(N = 10^{12}\), around 5–10 occurrences of sequences with 12 identical ending digits may be observed. Similarly, extending this reasoning further, one might predict approximately 5–10 occurrences of sequences with 15 consecutive identical ending digits up to \(N = 10^{15}\). The code used to compute these results is provided in \cite{supp}.

This numerical evidence does not prove the existence of arbitrarily long sequences of consecutive primes ending in the same digit. However, in 2000, Shiu proved the existence of sequences of any prescribed length \( k \) consisting of consecutive primes that all share the same remainder modulo \( q \)~\cite{shiu2000strings}. That is, there exist arbitrarily long sequences of the form:
\[
p_{n+1} \equiv p_{n+2} \equiv \dots \equiv p_{n+k} \equiv a \pmod{q},
\]
where \( a \) is an admissible residue modulo \( q \). In particular, for \( q = 10 \), the possible values of \( a \) are 1, 3, 7, or 9; corresponding to the only valid final digits of prime numbers in base 10.

This guarantees that sequences of consecutive primes ending in the same digit can be arbitrarily long. Consequently, the area \( A(N) \) is unbounded, since arbitrarily large steps in a fixed direction (e.g., always right, always up) are guaranteed to occur at some point in the walk, regardless of its current position.



It is worth mentioning that even though the area covered by the PWs increases indefinitely (a fact now established as a theorem rather than a conjecture) this does not guarantee continuous linear growth indefinitely. One might reasonably anticipate intervals in which this growth slows, potentially shifting from linear to sub-linear or even asymptotic behavior. Predicting precisely when and how such transitions might occur remains an intriguing open question.

\section{3D Structure and Benford's Law}

The coordinate \( z(x, y) \) was defined in~\cite{fraile2021primes} at each point \((x, y)\) as the total number of steps that the PW spent at that location, including all stationary steps between consecutive prime numbers. Under this definition, the \( z(x, y) \) values represented a mix of sums of visits and gaps between consecutive primes. Here, we adopt a different approach: we define \( z(x, y) \) simply as the number of visits to each point, explicitly excluding the stationary steps. Consequently, the sum of all \( z(x, y) \) values now equals the total number of primes \( N_p \).

A natural question arises regarding the statistical distribution of the \( z \)-values under our new definition: do they exhibit a random distribution, or do they follow a specific statistical pattern? Intriguingly, similar to the observations reported in our earlier study~\cite{fraile2021primes}, we find that these \( z_{\text{max}} \)-values clearly adhere to Benford’s Law~\cite{szpiro37}. Benford’s Law, one of the most celebrated laws governing leading-digit distributions~\cite{Hill1995, Hill1996, luque09}, states that the leading digit \( d \) (the leftmost nonzero digit of a number) is not uniformly distributed. Instead, it follows the logarithmic probability distribution
\[
P(d) = \log_{10}\left(1 + \frac{1}{d}\right), \quad d = 1, 2, \dots, 9.
\]

\begin{figure}
    \centering
    \includegraphics[width=0.5\textwidth]{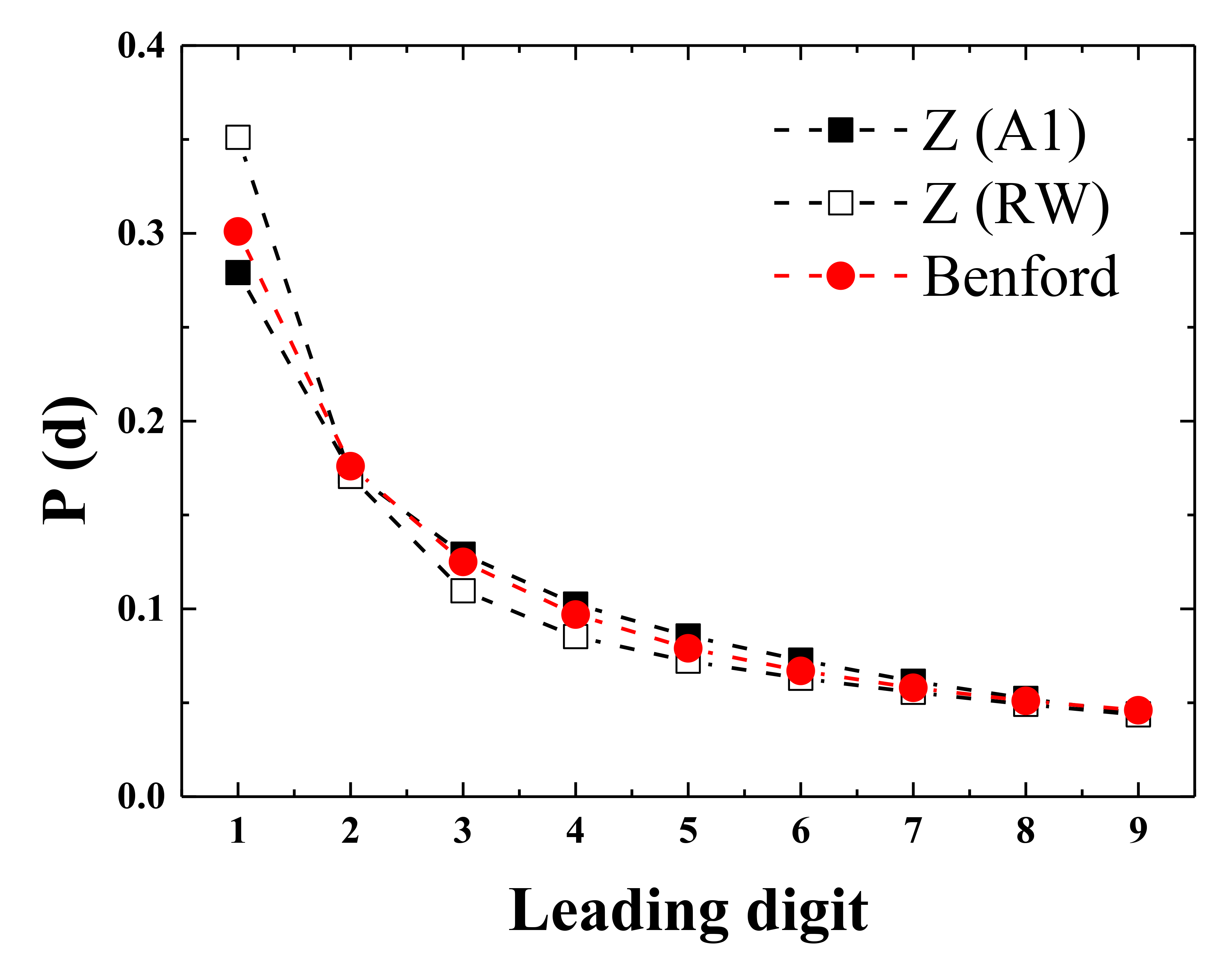}
    \caption{Leading-digit histogram of \( z_{\text{max}} \)-values for the PW up to \( N = 10^{10} \). Black squares indicate observed proportions for each leading digit, and red symbols show the expected Benford distribution. Open squares correspond to the results computed from a pure RW (number of steps \( = 10^9 \)). \( P(d) \): proportion; \( d \): leading digit.}
    \label{fig:benford_histogram}
\end{figure}

Figure~\ref{fig:benford_histogram} clearly illustrates this effect by plotting the leading-digit histogram of the \( z_{\text{max}} \)-values obtained for the PW up to \( N = 10^{10} \). The black squares depict the observed proportions of each leading digit, while the red symbols represent the expected values according to Benford’s Law. Even for a purely RW, which tends to produce a less compact region, the \( z \)-values also follow Benford’s Law closely (open squares). Further analysis of the three-dimensional structure of the prime walks is presented in the SM. For instance, the histograms of \( z_{\max} \) for each PW show near-overlapping curves, highlighting their shared statistical behavior despite differences in trajectory.

\section{Polar Representation of Prime Walks}

Exploring alternative representations of prime walks may reveal unique characteristics, particularly if there exist subtle statistical irregularities or biases within the prime digit sequences modulo 10. One might conceptualize an underlying hypothetical mapping responsible for generating such sequences. Previous studies have employed entropic and time-series analyses to study prime sequences~\cite{bonanno17}, while polar coordinates have been used effectively to analyze Brownian motion~\cite{godreche21}.

In our analysis, it is insightful to examine discrete differences of the radius, providing bounded oscillatory dynamics that can be directly compared to angular changes per step of the walk. We define the radius and angle as
\[
R_n = \sqrt{X_n^2 + Y_n^2} \quad \text{and} \quad \phi_n = \tan^{-1}\left(\frac{Y_n}{X_n}\right) \, .
\]

\begin{figure}[h]
    \centering
    \includegraphics[width=0.5\textwidth]{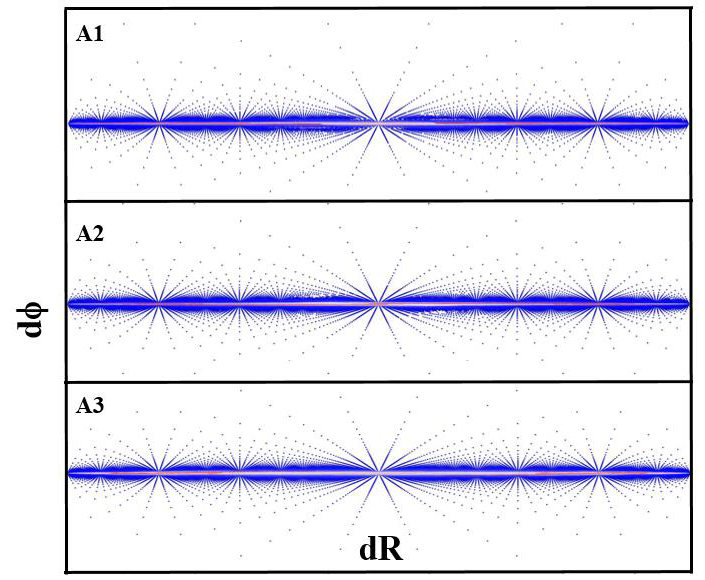}
    \caption{Polar plot of $\Delta\phi$ versus $\Delta R$ for the three prime walks (A1: top, A2: middle, A3: bottom), computed up to $N = 10^6$. The color gradient represents the step number to enhance visual clarity.}
    \label{fig:fractals1}
\end{figure}

Initially, comparing direct discrete differences of the radius, \( \Delta R = R_{n+1} - R_n \), to angular positions resulted in low-frequency cyclic patterns reminiscent of Lissajous figures, complicating interpretation. To address this, we instead compared discrete differences of both quantities: \( \Delta R \) versus \( \Delta \phi = \phi_{n+1} - \phi_n \). This approach produced a more linear mapping, revealing unexpected fractal characteristics in the resulting point-cloud distribution (see Figure~\ref{fig:fractals1}).

A comparative analysis was conducted between our original PWs and an unbiased, four-step planar Pearson-type RW, wherein each step along the \( X \) and \( Y \) axes was selected with equal probability from a uniform random number generator, incorporating a discrete selection criterion of the form \(\lfloor r_n/0.25 \rfloor\).

Figure~\ref{fig:fractals1} display the point clouds generated by the three variants of PWs. Direct comparison indicates that the RW and the third type of PW exhibit relatively suppressed fractal characteristics, whereas the first two PWs present more pronounced fractal-like point clouds. This observation suggests potentially greater structural complexity or statistical discrepancies inherent to these two PW algorithms.


A histogram of $\Delta\phi$ values for the three PWs and a RW reveals no significant differences between the distributions (see \cite{supp}).


\section{The Prime Walk Recurrence Conjecture}

The prime walker moves erratically; however, despite the continuously increasing explored area \( A(N) \), it never seems to venture too far from the origin. Figure~\ref{fig:area} clearly illustrates the most visited regions for each of the three PW definitions up to \( N = 10^{10} \). Notably, for algorithms A1 and A2, the most frequently visited point is located relatively close to the origin, and for algorithm A3, it is situated even closer. This simple yet insightful observation motivates the following straightforward conjecture:

\textit{In the limit as \( N \to \infty \), the most frequently visited point \((x, y)\) by the PW is the origin \((0, 0)\).}

This likely holds true for all three PW algorithms - see Figure~2 in \cite{supp}, where two additional conjectures are presented.

One concerns the asymptotic distribution of the area explored by the PWs, suggesting that it becomes evenly spread across all four quadrants. The other proposes that the distance from the origin to the most frequently visited points, as well as to those with maximal $z$-values, remains bounded by finite constants. These conjectures are motivated by numerical observations and highlight potential structural regularities in the behavior of PWs. 

\section{Conclusions}
\label{sec:Conclusion}

The results presented in this work demonstrate that the area traversed by the PW is unbounded. This is supported by the existence of arbitrarily long sequences of consecutive primes ending in any of the four allowed terminal digits \cite{shiu2000strings}, ensuring that the PW continues its growth indefinitely.

We have shown that the three PW algorithms, despite differing in their directional assignment rules, yield equivalent behavior. Specifically, with an explored area that grows linearly with the same slope, \(\beta = 0.03\), across all cases.  We interpret this consistency as strongly supportive of the hypothesis that linear growth persists regardless of the interval length considered.

Our study notably confirms complete alignment with Benford’s Law, despite the \( z \)-values in this analysis being derived through a significantly different approach compared to previous research. This robustness enhances the credibility of the statistical patterns observed. 

Beyond extending the PW framework established in ref.~\cite{fraile2021primes}, our findings open several intriguing research directions. First, they allow a fresh investigation of the classical question concerning runs of consecutive primes that share the same terminal digit. Second, they raise deeper questions about the underlying symmetries and statistical structures embedded in the distribution of prime numbers.

The framework also encourages a broader investigation into prime-based random walks. In ref.~\cite{fraile2020jacobs}, Jacob’s Ladder is introduced,
a straightforward one-dimensional walk derived from primes, named for its alternating pattern of ascent and descent. Despite its apparent simplicity, it raises profound mathematical questions, most notably whether it crosses zero infinitely often: a problem that may remain unsolved.

These investigations emphasize a broader philosophical insight: experimental mathematics is not a peripheral activity but a core part of the mathematical research method. With modern computational tools, we are now equipped to explore the primes and their hidden structures with unprecedented depth and scale \cite{borwein2003experiment, borwein2004experimentation}.


\section*{Credit authorship contribution statement}
\label{sec:CRediT}

\noindent \textbf{A. Fraile:} Writing - original draft, investigation, formal analysis, data collection, conceptualization, funding acquisition. \\
\textbf{D. Fernández:} Writing – review \& editing, validation. \\
\textbf{R. Martínez:} Coding, review \& editing, validation. \\
\textbf{T. Raptis:} Coding, review \& editing, validation.

\section*{Declaration of competing interest}
\label{sec:CompetingInterest}

\noindent The authors declare that they have no competing financial interests or personal relationships that could have influenced the work reported in this paper.

\section*{Acknowledgments}
\label{sec:Acknowledgments}

\noindent The authors are indebted to Prof. James Tuite for several insightful discussions. AF acknowledges support from ICN2, which is supported by the Severo Ochoa programme from Spanish MINECO (grant CEX2021-001214-S) and by Generalitat de Catalunya (CERCA programme and Grant 2021SGR01519). DF acknowledges support from Rannís (Icelandic Centre for Research) through a Sproti grant (Project No. 2423212-601).

\bibliography{references}

\end{document}